\documentclass[a4paper,12pt]{article}
\usepackage[english]{babel}
\usepackage{graphicx}
\usepackage{amssymb}
\usepackage{amsmath}
\usepackage{geometry}
\usepackage{hyperref}
\hypersetup{
colorlinks=true,
linkcolor=black,
citecolor=black,
urlcolor=black
}

\begin{document}

\newcommand{\e}{{\sf E}}
\newcommand{\D}{{\sf D}}
\newcommand{\p}{{\sf P}}
\newcommand{\N}{\mathbb N}
\renewcommand{\r}{\mathbb R}
\renewcommand*{\Pr}{\mathbb P}
\newcommand*{\R}{\mathbb R}
\newcommand*{\I}{\mathbb I}

\renewcommand{\refname}{References}

\title{\vspace{-1.8cm}
Statistical Analysis of Precipitation Events}

\author{
V.\,Yu.~Korolev\textsuperscript{1},
A.\,K.~Gorshenin\textsuperscript{2},
S.\,K.~Gulev\textsuperscript{3},
K.\,P.~Belyaev\textsuperscript{4},
A.\,A.~Grusho\textsuperscript{5}}

\date{}

\maketitle

\footnotetext[1]{Faculty of Computational Mathematics and
Cybernetics, Lomonosov Moscow State University, Russia; Institute of
Informatics Problems, Federal Research Center ``Computer Science and
Control'' of Russian Academy of Sciences, Russia; Hangzhou Dianzi
University, China; \url{vkorolev@cs.msu.su}}

\footnotetext[2]{Institute of Informatics Problems, Federal Research
Center ''Computer Science and Control'' of Russian Academy of
Sciences, Russia; \url{agorshenin@frccsc.ru}}

\footnotetext[3]{P.P. Shirshov Institute of Oceanology, Russia; \url{gul@sail.msk.ru}}

\footnotetext[4]{P.P. Shirshov Institute of Oceanology, Russia; \url{kosbel55@gmail.com}}

\footnotetext[5]{{Institute of Informatics Problems, Federal Research Center ``Computer Science and Control'' of Russian Academy of Sciences, Russia; Faculty of Computational Mathematics and Cybernetics, Lomonosov Moscow State University, Russia; \url{ grusho@yandex.ru}}}

\maketitle

\begin{abstract}
In the present paper we present the results of a statistical
analysis of some characteristics of precipitation events and propose
a kind of a theoretical explanation of the proposed models in terms
of mixed Poisson and mixed exponential distributions based on the
information-theoretical entropy reasoning. The proposed models can
be also treated as the result of following the popular Bayesian
approach.
\end{abstract}

\section{Introduction}

In most papers available to the authors, in which meteorological
data is analyzed statistically, the suggested analytical models for
the observed statistical regularities in precipitation are rather
ideal and far from being adequate. For example, it is traditionally
assumed that the duration of a wet period (the number of
subsequent wet days) follows the geometric distribution (for example, see~\cite{Zolina2013}). Perhaps,
this prejudice is based on the conventional interpretation of the
geometric distribution in terms of the Bernoulli trials as the
distribution of the number of subsequent wet days
(``successes'') till the first dry day (``failure''). But the
framework of Bernoulli trials assumes that the trials are
independent whereas a thorough statistical analysis of precipitation
data registered in different points demonstrates that the sequence
of dry and wet days is not only independent, but it is also
devoid of the Markov property so that the framework of Bernoulli
trials is absolutely inadequate for analyzing meteorological data.

In the present paper we present the results of a statistical
analysis of some characteristics of precipitation events and propose
a kind of a theoretical explanation of the proposed models in terms
of mixed Poisson and mixed exponential distributions based on the
information-theoretical entropy reasoning. The proposed models can
be also treated as the result of following the popular Bayesian
approach. The adequate models of the regularities in precipitation
events are very important for adequate forecasting of natural
disasters such as water floods or long dry spells.

\section{The analysis of the duration of wet periods}

\begin{figure}[h]
\begin{minipage}[h]{0.49\textwidth}
\center{\includegraphics[width=\textwidth]{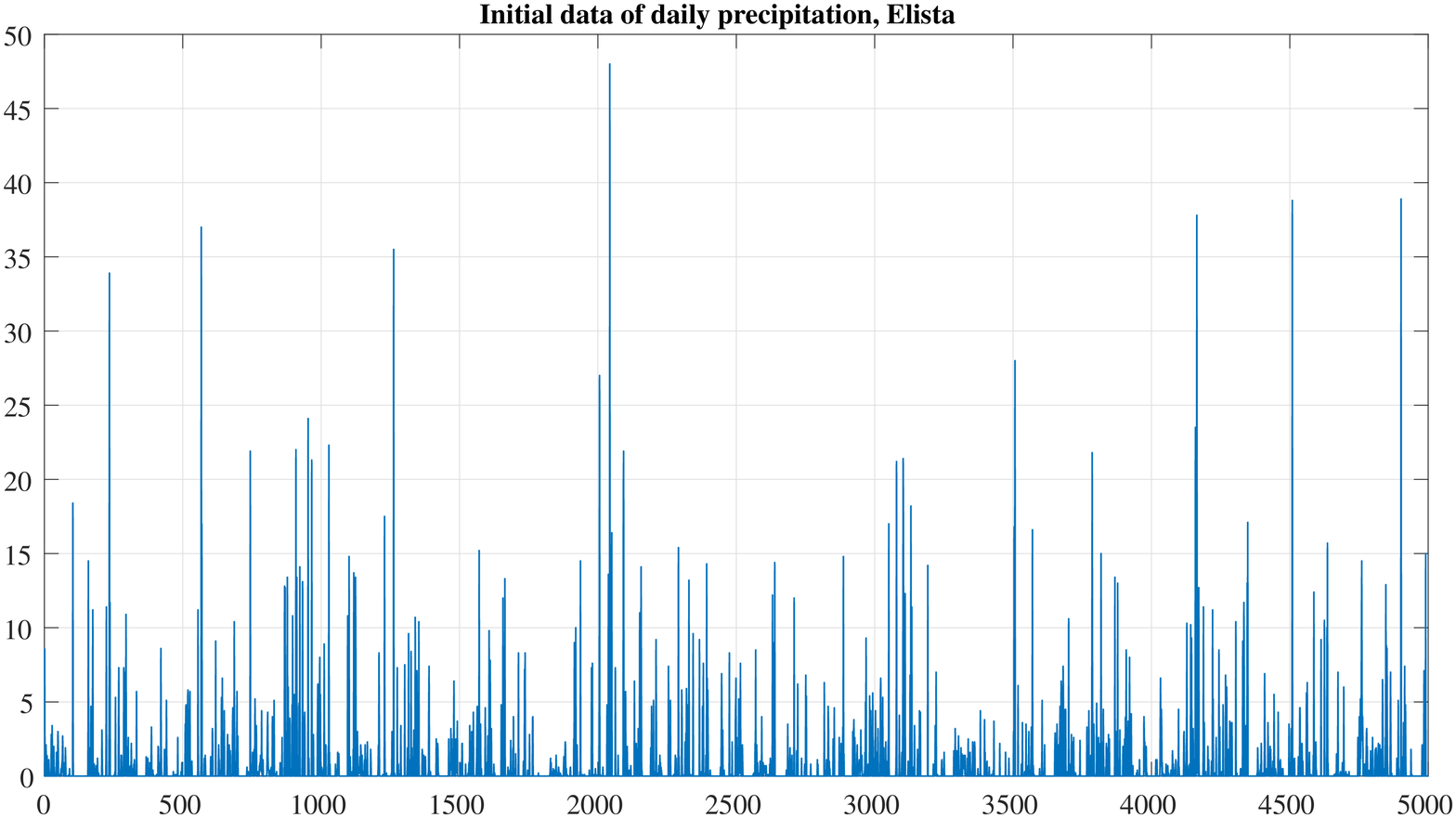} \\a)}
\end{minipage}
\hfill
\begin{minipage}[h]{0.49\textwidth}
\center{\includegraphics[width=\textwidth]{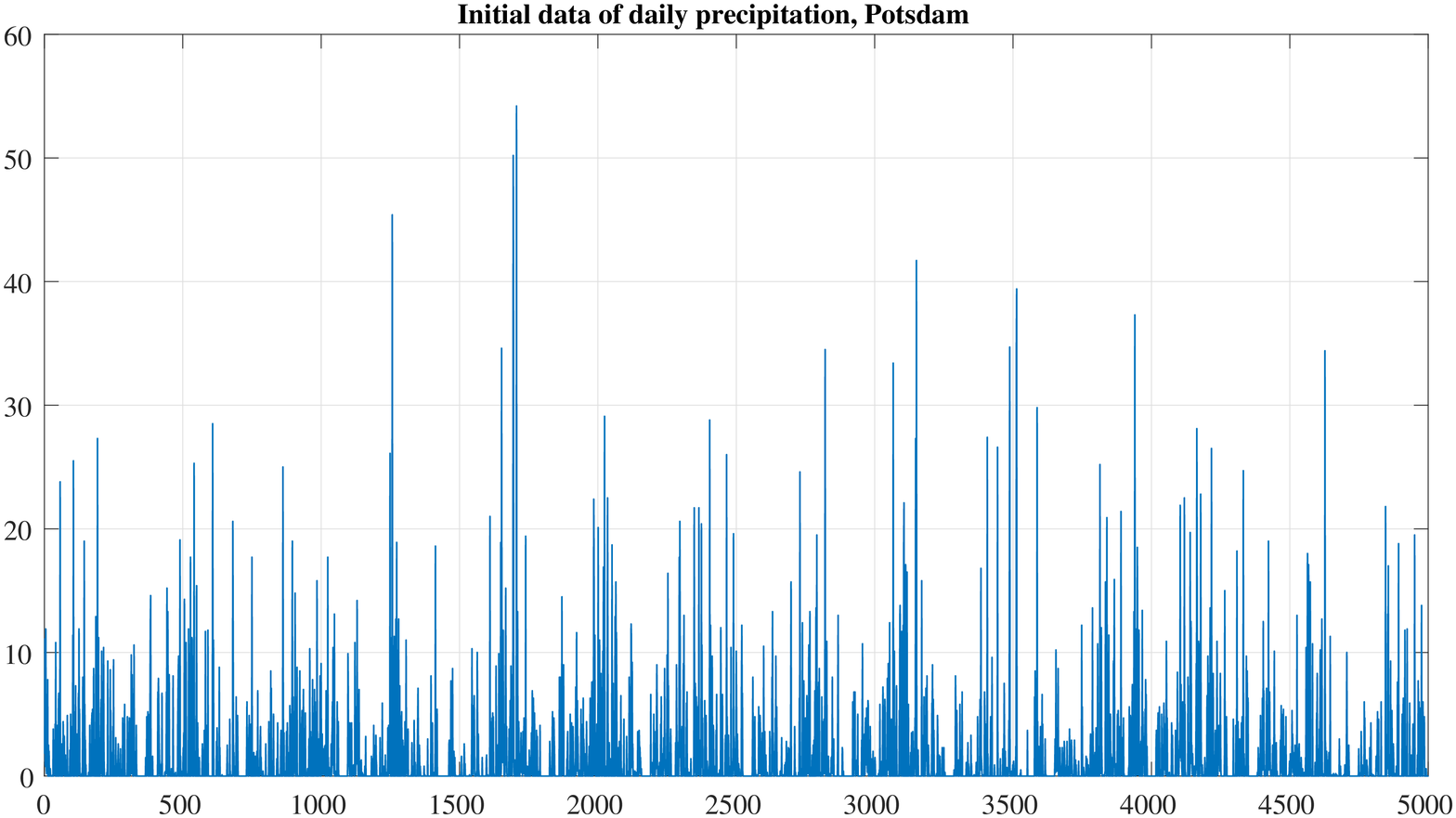} \\ b)}
\end{minipage}
\label{Data}
\caption{The initial data of daily precipitation in Elista (a) and Potsdam (b).}
\end{figure}

We analyzed meteorological data registered at two point with very
different climate: Potsdam (Germany) with mild climate influenced by
the closeness to the ocean with warm Gulfstream flow and Elista
(Russia) with radically continental climate. The initial data of daily precipitation in Elista and Potsdam are
presented on Figure~\ref{Data}a and Figure~\ref{Data}b, respectively.
On these figures the horizontal axis is discrete time
measured in days. The vertical axis is the daily precipitation
volume measured in centimeters. In other words, the height of each
``pin'' on these figures is the precipitation volume registered at
the corresponding day (at the corresponding point on the horizontal
axis).

\begin{figure}[h]
\begin{minipage}[h]{0.49\textwidth}
\center{\includegraphics[width=\textwidth]{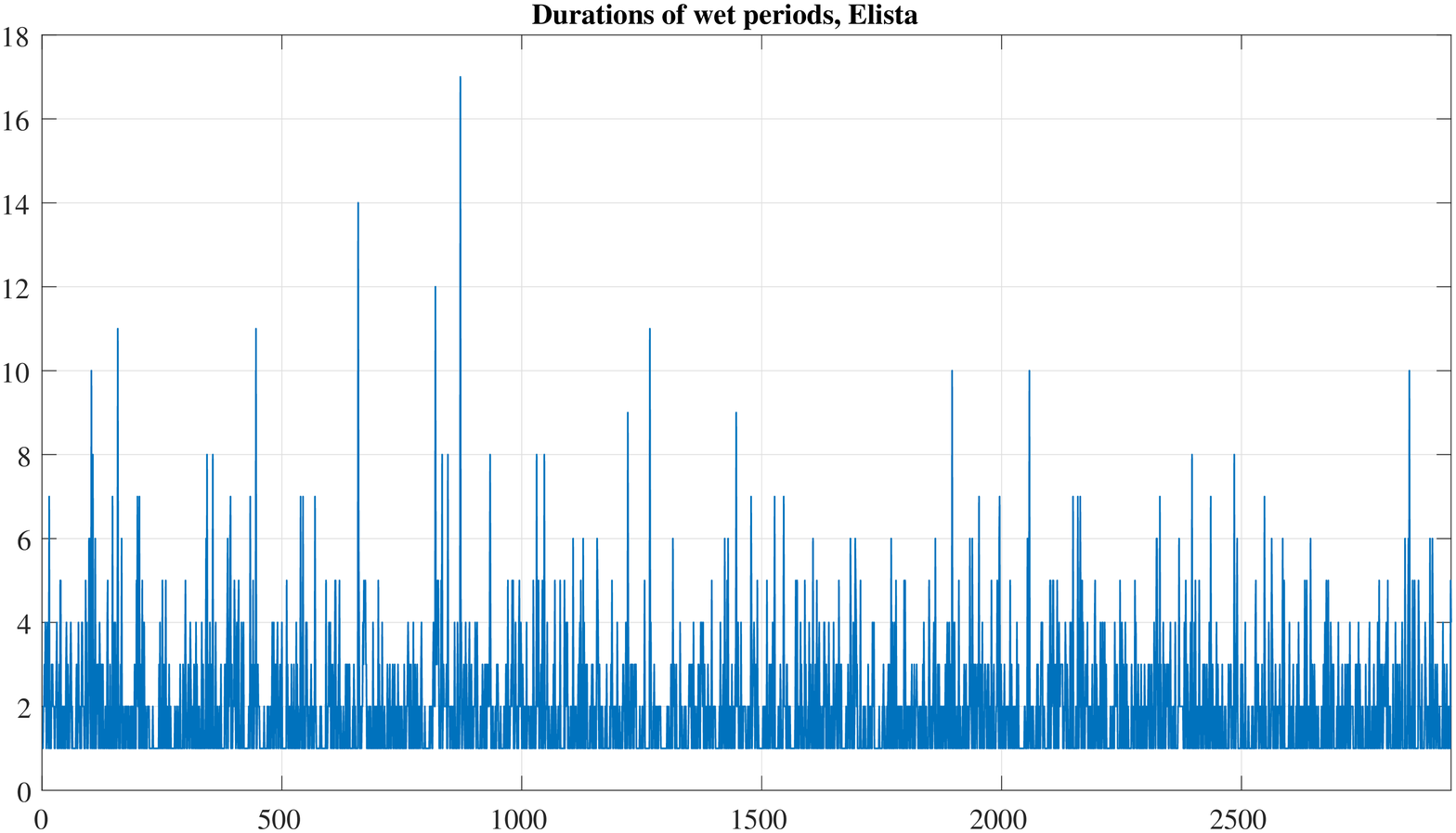} \\a)}
\end{minipage}
\hfill
\begin{minipage}[h]{0.49\textwidth}
\center{\includegraphics[width=\textwidth]{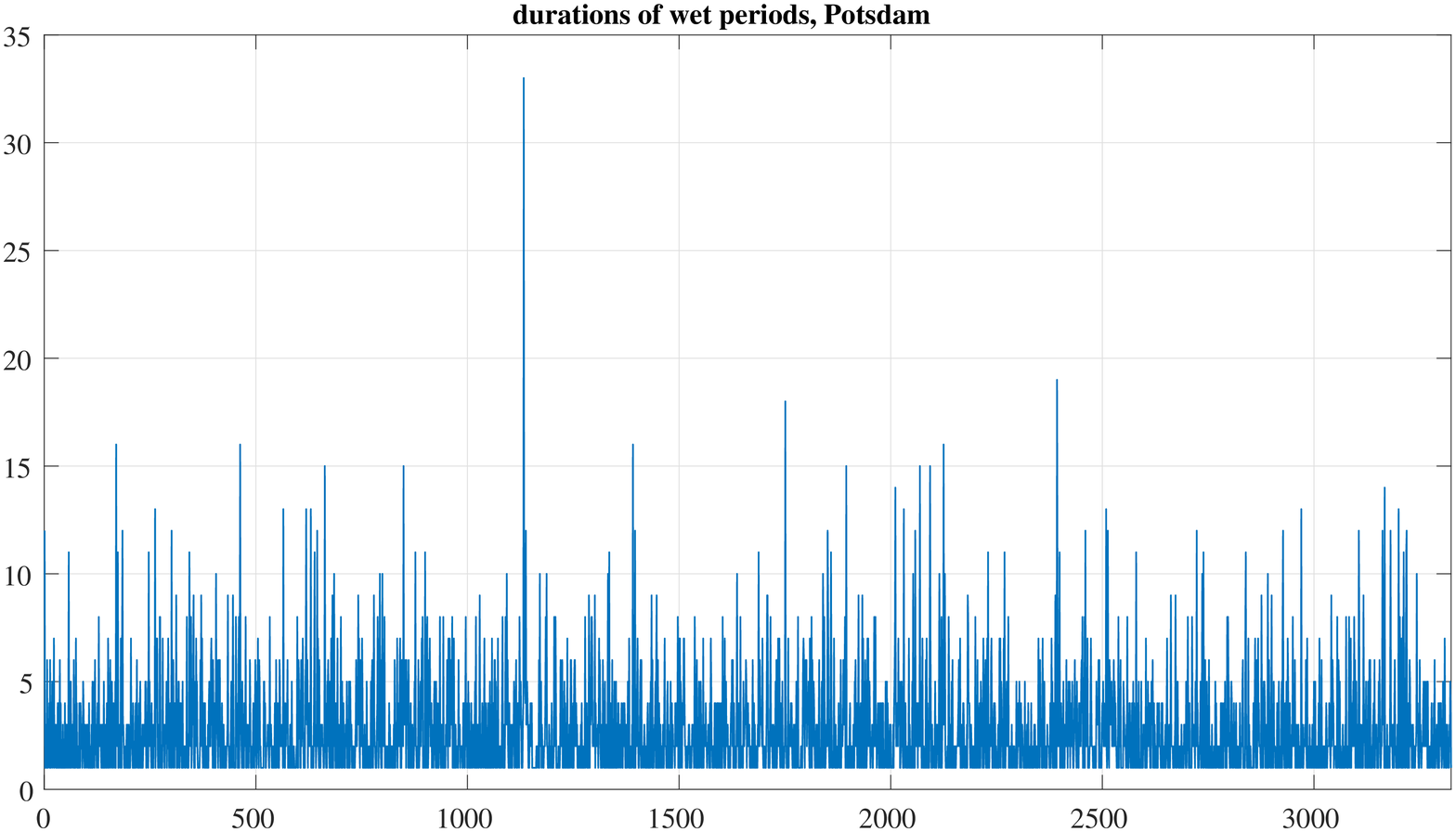} \\ b)}
\end{minipage}
\label{WetPeriod}
\caption{The durations of wet periods in Elista (a) and Potsdam (b).}
\end{figure}

In order to analyze the statistical regularities of the duration of
wet periods this data was rearranged as shown on Figure~\ref{WetPeriod}a and Figure~\ref{WetPeriod}b.

On these figures the horizontal axis is the number of successive
wet periods. It should be mentioned that directly before and after
each wet period there is at least one dry day, that is,
successive wet periods are separated by dry periods. On the vertical
axis there lie the durations of wet periods. In other words, the
height of each ``pin'' on these figures is the length of the
corresponding wet period measured in days and the corresponding
point on the horizontal axis is the number of the wet period.

The samples of durations in both Elista and Potsdam were assumed
homogeneous and independent. The best fit with these samples was
demonstrated by the negative binomial distribution. Let $r>0$ and
$p\in(0,1)$. We say that a random variable $Y$ has the negative
binomial distribution with parameters $r$ and $p$, $Y\sim NB(r,p)$,
if
\[
\p(Y=k)=\frac{\Gamma(r+k)}{\Gamma(r)\Gamma(k+1)} p^r (1-p)^{k},
\quad k=0,1,2,\ldots
\]
Let $X=Y+1$. Then $X\geqslant 1$ and
\[
\p(X=k)=\p(X-1=k-1)=\frac{\Gamma(r+k-1)}{\Gamma(r)\Gamma(k-1+1)} p^r
(1-p)^{k-1}=\frac{\Gamma(r+k-1)}{\Gamma(r)\Gamma(k)} p^r
(1-p)^{k-1}, 
\]
$k=1,2,3,\ldots$

\begin{figure}[h]
\begin{minipage}[h]{0.49\textwidth}
\center{\includegraphics[width=\textwidth]{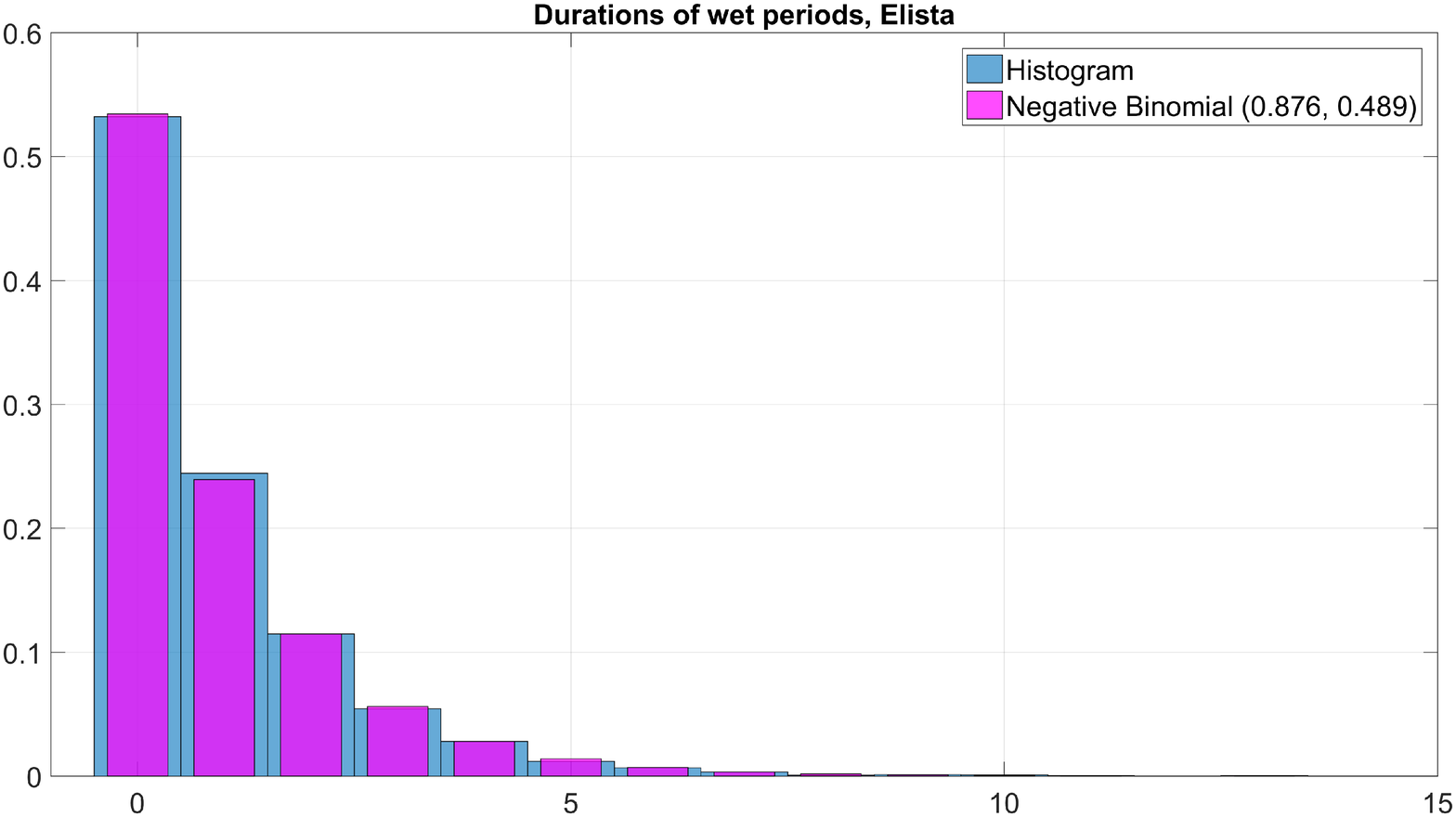} \\a)}
\end{minipage}
\hfill
\begin{minipage}[h]{0.49\textwidth}
\center{\includegraphics[width=\textwidth]{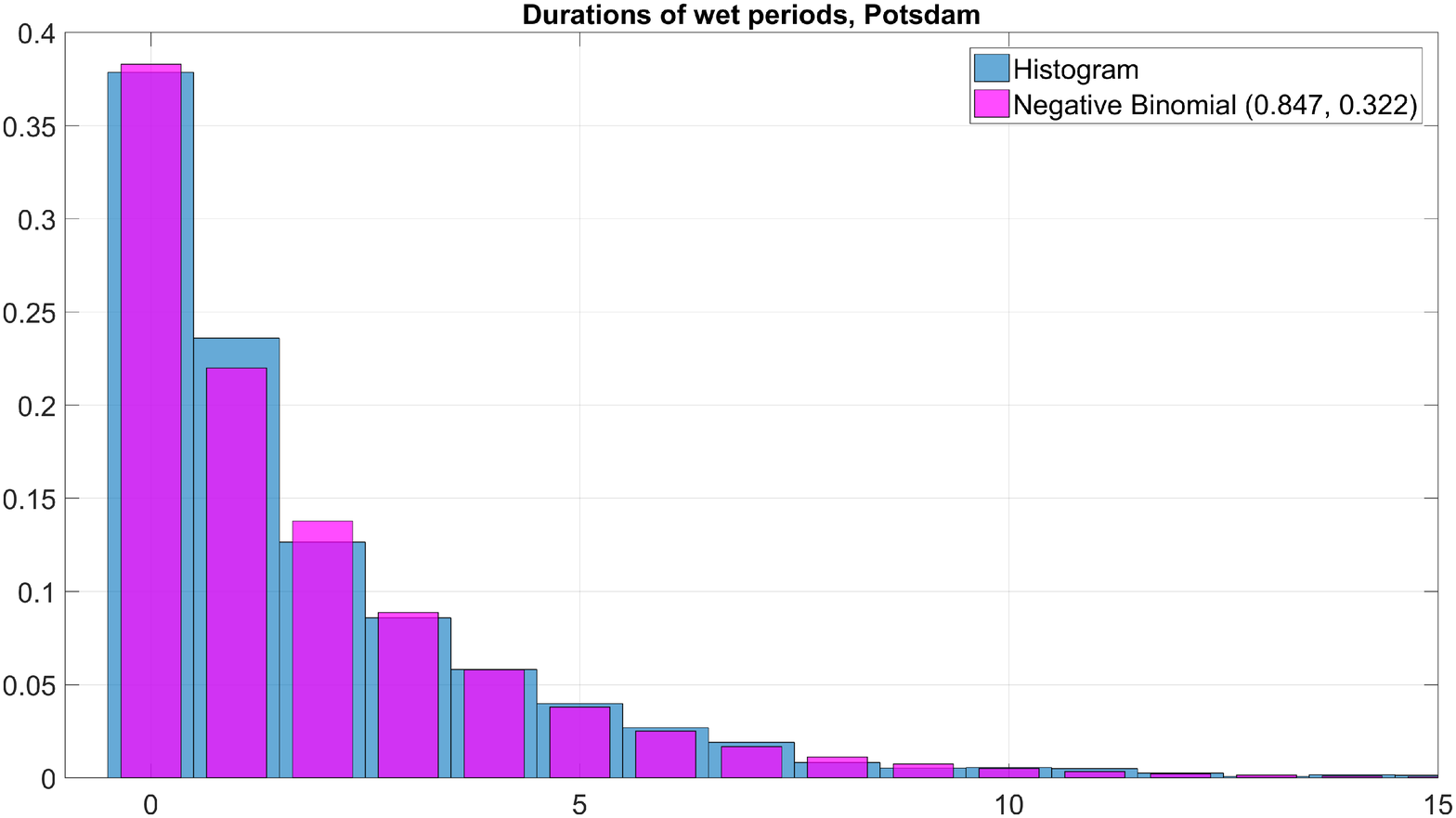} \\ b)}
\end{minipage}
\label{WetHist}
\caption{The histogram of durations of wet periods in Elista (a) and Potsdam (b) and the
fitted negative binomial distribution.}
\end{figure}

Figures~\ref{WetHist}a and~\ref{WetHist}b show the histograms
constructed from the corresponding samples of duration periods and
the fitted negative binomial distribution. In both cases the shape
parameter $r$ turned out to be less than one. For Elista $r=0.876$,
$p=0.489$, for Potsdam $r=0.847$, $p=0.322$.

We can suggest the following explanation for the negative binomial
distribution of the duration of wet periods. Each precipitation
process develops in temporal and spatial coordinates. Assume that
the projection of each precipitation event onto the time axis is a
(random) connected set. Then attracting some rather simple, natural
and well-tractable assumptions defining Poisson random measures (for example, see~\cite{Kingman1993}), we come to the suggestion that, if the
``temporal density'' of wet days is uniform in time, then the size
of each temporal projection of a precipitation event must have the
Poisson distribution. But due to seasonal and other climatic trends
the ``temporal density'' of wet days varies in a random manner.
Therefore, the actual distribution of the duration of wet periods
should be sought among mixed Poisson distributions. As far ago as in
1920 it was shown that the negative binomial distribution is mixed
Poisson with the mixing gamma-distribution~\cite{GreenwoodYule1920}:
for $k=0,1,\ldots$ we have
\begin{gather}
\frac{1}{k!}\int_{0}^{\infty}e^{-\lambda}\lambda^k\frac{1}{\Gamma(r)}\Big(\frac{p}{1-p}\Big)^r\lambda^{r-1}\exp\Big\{-\frac{\lambda
p}{1-p}\Big\}d\lambda=\notag\\
\Big(\frac{p}{1-p}\Big)^r\frac{1}{k!\Gamma(r)}\int_{0}^{\infty}\exp\Big\{-\frac{\lambda}{1-p}\Big\}\lambda^{k+r-1}d\lambda=
\frac{\Gamma(k+r)}{k!\Gamma(r)}p^r(1-p)^k. \label{NBmixedPosson}
\end{gather}
In other words, if the parameter $\lambda$ of the Poisson
distribution is random and has the gamma-distribution with shape
parameter $r$ and scale parameter $p/(1-p)$, then the resulting
mixed Poisson distribution is negative binomial with parameters $r$
and $p$.

Note that in both cases, that is, for Elista and Potsdam, the
parameter $r$ of the fitted negative binomial distribution and,
hence, of the mixing gamma-distribution in (\ref{NBmixedPosson}), is
less than one. It is known that the gamma-distribution can be
represented as a mixture of exponential laws if and only if its
shape parameter is no greater than one (see~\cite{Gleser1989}). In
that paper it was proved that if the shape parameter $r$ of the
gamma-distribution $f_{r,\theta}(x)$ satisfies the condition $0<r\leq 1$,
then the gamma-density $f_{r,\theta}(x)$ can be represented as a
mixed exponential distribution:
\begin{equation}
\label{CondForm} f_{r,\theta}(x)=\int_{0}^{\infty}p_\theta(\gamma)
\gamma e^{-\gamma x}\,d\gamma.
\end{equation}
where
\begin{equation}\label{pForm}
p_\theta(\gamma)
=\frac{(\gamma-\theta)^{-r}\theta^r}{\gamma\Gamma(1-r)\Gamma(r)}
\I(\gamma\geqslant\theta).
\end{equation}
But it is well known that the exponential distribution has the
maximum differential entropy among all distributions concentrated on
the nonnegative half-line and possessing a finite expectation. This
means that the negative binomial model for the duration of wet
periods has serious theoretical grounds leading to the reasoning
based on the principle of non-decrease of uncertainty in closed
systems. But any regional precipitation system cannot be considered
as closed, it is influenced by many poorly predictable random
factors. Within the mixed models (\ref{NBmixedPosson}) and
(\ref{CondForm}) the distribution $p_\theta(\gamma)$~(\ref{pForm}) accumulates the
information concerning the statistical regularities in the integral
behavior of ``external'' factors.

\section{The analysis of the daily precipitation volumes}

\begin{figure}[h]
\begin{minipage}[h]{0.49\textwidth}
\center{\includegraphics[width=\textwidth]{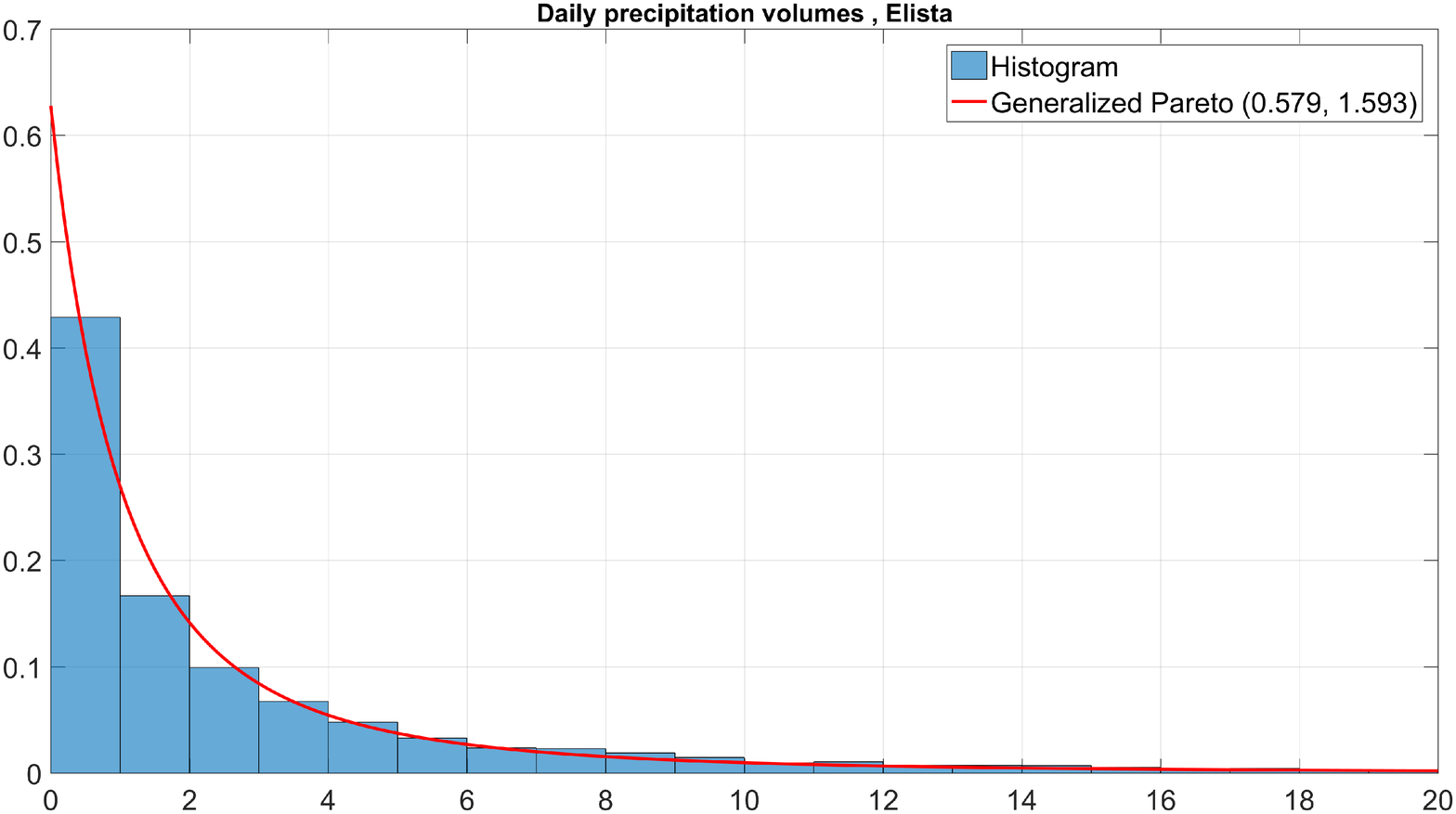} \\a)}
\end{minipage}
\hfill
\begin{minipage}[h]{0.49\textwidth}
\center{\includegraphics[width=\textwidth]{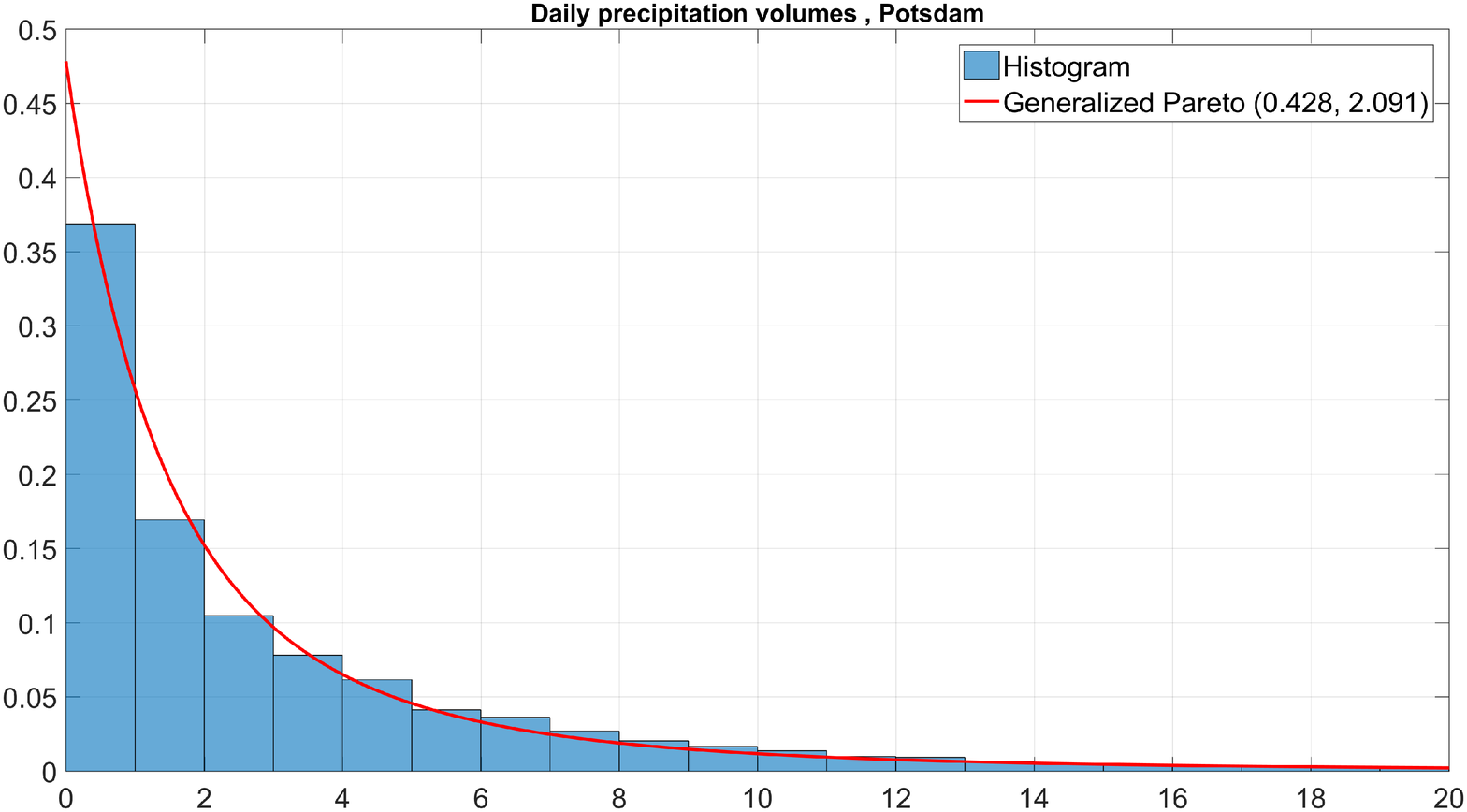} \\ b)}
\end{minipage}
\label{Precip}
\caption{The histogram of daily precipitation volumes in Elista (a) and Potsdam (b) and the fitted Pareto distribution.}
\end{figure}

The distribution of daily precipitation volumes presented on
Figure~\ref{Precip}a and Figure~\ref{Precip}b was also analyzed. It
turned out that for both Elista and Potsdam data the best fit is
demonstrated by the Pareto distribution defined by the distribution
function
\[
F_{\xi, \sigma, \mu}(x)=
1-\left(1+\frac{\xi(x-\mu)}{\sigma}\right)^{-1/\xi}, \quad \xi\neq
0,
\]
see Figure~\ref{Precip}a and Figure~\ref{Precip}b.

As concerns the Pareto distribution, it can be easily seen that
scale mixtures of gamma-distributions with the mixing law also
having the gamma-distribution are Pareto-type distributions. Indeed,
we have
\[
\int_{0}^{\infty}\lambda e^{-\lambda
x}\frac{\mu^s}{\Gamma(s)}\lambda^{s-1}e^{-\mu\lambda}d\lambda=
\frac{s\mu^s}{(x+\mu)^{1+s}}, \ \ \ x>0,
\]
that is, the Pareto distribution with the density
\[
p(x;s,\mu)=s\mu^s(x+\mu)^{-(1+s)},\ \ \ \ x>0,
\]
for any $s>0$, $\mu>0$ is a mixed exponential distribution.

Moreover, a more general result can be easily proved: all scale
mixtures of gamma-distributions with the mixing law also having the
gamma-distribution are Pareto-type distributions. Indeed, we have
\[
\frac{x^{r-1}}{\Gamma(r)}\int_{0}^{\infty}\lambda^re^{-\lambda
x}\frac{\mu^s}{\Gamma(s)}\lambda^{s-1}e^{-\mu\lambda}d\lambda=
\frac{x^{r-1}\mu^s}{\Gamma(r)\Gamma(s)}\int_{0}^{\infty}\lambda^{r+s-1}e^{-\lambda(x+\mu)}d\lambda=
\]
\[
=\frac{x^{r-1}\mu^s}{\Gamma(r)\Gamma(s)(x+\mu)^{r+s}}\int_{0}^{\infty}(x+\mu)^{r+s-1}\lambda^{r+s-1}e^{-\lambda(x+\mu)}d\lambda(x+\mu)=
\frac{\Gamma(r+s)\mu^s}{\Gamma(r)\Gamma(s)}\frac{x^{r-1}}{(x+\mu)^{r+s}},\
\ \ x>0.
\]

And if here the shape parameter $s$ of the mixing gamma-distribution
satisfies the condition $s\leqslant1$, then from the abovesaid it follows
that such Pareto distributions are also mixed exponential so that
the reasoning presented in the final part of the preceding section
can be applied to the daily precipitation volumes.

We also propose some similar mixture-type models for other
characteristics of precipitation process, such as total
precipitation volume during a wet period, etc.

\section*{Acknowledgments}
The research is partially supported by the Russian Foundation for
Basic Research (projects 15-37-20851 and 15-07-04040) and the government project of FRC CSC RAS No. 0063-2015-0014.

\bibliographystyle{aipnum-cp}%

\end{document}